\documentclass{ifacconf}

\usepackage{graphicx}
\usepackage{natbib}

\usepackage{amsmath,amssymb}
\usepackage{booktabs}
\usepackage{tabularx}
\usepackage{array}

\usepackage{algorithm}
\usepackage{algpseudocode}

\usepackage{xcolor}

\newcommand*{\argmin}{\operatornamewithlimits{argmin}\limits}

\newcommand*{\argminOp}{\operatornamewithlimits{argmin}\limits}

\makeatletter
\let\old@ssect\@ssect 
\makeatother

\usepackage[bookmarks=false]{hyperref}
\hypersetup{
    colorlinks = true,
    linkcolor  = {magenta!90!black},
    citecolor  = {cyan!70!black},
    urlcolor   = {blue!20!black}
}
\allowdisplaybreaks
\edef\endfrontmatter{%
  \unexpanded\expandafter{\endfrontmatter}%
  \noexpand\endNoHyper
}
\makeatletter
\def\@ssect#1#2#3#4#5#6{%
  \NR@gettitle{#6}
  \old@ssect{#1}{#2}{#3}{#4}{#5}{#6}
}
\makeatother
\allowdisplaybreaks
\begin{document}
\begin{frontmatter}

\title{Stochastic Optimization for Switch-System Identification under Errors-in-Variables} 

\author[First]{Yasin Gharehmohammadlou}
\author[Second]{Mohamad Kazem Shirani Faradonbeh}
\author[Third]{Mohammad Khosravi}
\author[Last]{Digvijay Boob}

\address[First]{Operations Research and Engineering Management, Lyle School of Engineering,
Southern Methodist University,
Dallas, Texas, USA (e-mail: 	
ygarehmohammadlou@smu.edu)}
\address[Second]{Department of Mathematics,
Southern Methodist University,
Dallas, Texas, USA (e-mail: 	
mohamadksf@smu.edu)}
\address[Third]{Delft Center for Systems and Control, Department of Mechanical Engineering, Delft University of Technology, Delft, The Netherlands (e-mail: 
mohammad.khosravi@tudelft.nl)}
\address[Last]{Operations Research and Engineering Management, Lyle School of Engineering,
Southern Methodist University,
Dallas, Texas, USA (e-mail: 	
dboob@smu.edu)}

\begin{abstract}
 In this paper, we address the identification of switched dynamical systems under errors-in-variables, where both input and output measurements are corrupted by noise, and the number of operating modes and subsystem orders are unknown. We propose a unified identification framework that simultaneously estimates the subsystem parameters, latent switching sequence, and model structure. The proposed formulation combines a variance-normalized data-fitting criterion with temporally regularized soft mode assignments and sparsity-inducing penalties for selecting the active modes and relevant input-output lags. To solve the resulting nonconvex optimization problem, we develop an alternating algorithm that incorporates convex mode-assignment updates and proximal stochastic variance-reduced gradient optimization for parameter estimation. While the existing approaches in the literature require exhaustive searches over candidate structures, the proposed method performs structure selection within a single model. Numerical experiments demonstrate accurate recovery of the underlying system structure and favorable computational efficiency compared with a generalized total least squares and dynamic programming approach, particularly as the number of candidate structures increases.
\end{abstract}

\end{frontmatter}

\section{Introduction}
\label{sec:introduction}
Identification of switched dynamical systems with inputs
requires estimating subsystem parameters and unobserved
switching modes. These systems provide a flexible framework
for representing processes whose dynamics vary between
different operating regimes. The identification problem
involves the simultaneous estimation of continuous
subsystem parameters and discrete mode assignments,
introducing considerable computational challenges.
Comprehensive reviews of hybrid and switched-system
identification methods are provided in
\citet{paoletti2007tutorial,garulli2012survey},
while \citet{lauer2019hybrid} discusses the theoretical
and algorithmic foundations of learning switching models.

Early algebraic methods address switched-system
identification through polynomial representations
\citep{vidal2003algebraic,mavidal2005varieties},
while mixed-integer, sparse, and sum-of-norms
approaches provide optimization-based alternatives
\citep{roll2004pwa,bako2011sparse,ohlsson2013sum}.
An alternating framework for fitting jump models
through continuous parameter estimation and discrete
mode assignment is developed in
\citet{bemporad2018jump}.
Similarly, \citet{bianchi2020alternating} proposes
an alternating optimization approach that incorporates
prior information about the switching mechanism.
These approaches demonstrate the effectiveness of
optimization-based identification, although the
coupling between subsystem estimation and mode
assignment remains a central computational challenge.

Structural uncertainty introduces further difficulties
when the number of operating modes and the subsystem
orders are unknown. Randomized model-selection
approaches for switched systems are investigated in
\citet{bianchi2021switchnarx}, while statistical
learning methods for estimating the number of modes
are considered in \citet{massucci2022statistical}.
In addition, \citet{mu2022persistence} studies
persistence-of-excitation conditions that guarantee
uniqueness of the identified subsystem parameters.
More recently, \citet{zheng2024twostage} proposes
a two-stage identification framework combining
dynamic programming for switching-time estimation
and sparsity-inducing optimization for subsystem
identification. These studies highlight the importance
of structural modeling and computational efficiency.
Nevertheless, jointly determining the number of modes
and the input-output lag orders can require expensive
structural search, particularly when the admissible
upper bounds are large.

The presence of measurement noise further complicates
system identification. When both inputs and outputs
are contaminated by measurement errors,
errors-in-variables (EIV) estimation becomes necessary,
since ordinary least squares estimators generally
produce biased estimates when the regressors are noisy
\citep{markovsky2007tls}. Scalable identification
methods for switched-system models using large noisy
datasets are investigated in
\citet{hojjatinia2020noisy}, while
\citet{sznaier2024riemannian} proposes a Riemannian
embedding approach to EIV switched-system
identification. These contributions illustrate
the importance of noise-aware estimation, although
their modeling assumptions and identification
procedures differ. The simultaneous treatment of
measurement uncertainty, latent switching sequences,
and unknown system structure remains computationally
challenging.

Related advances in regularized system identification
include the kernel-based formulations in
\citet{khosravi2023uniqueness,khosravi2023frequency},
which investigate well-posedness and the incorporation
of prior system information. An expectation-maximization
approach for identifying linear dynamics with bilinear
observation models is developed in
\citet{liu2024bilinear}. From a statistical
perspective, \citet{faradonbeh2018finite} establishes
finite-time parameter-identification bounds for
unstable linear systems, while
\citet{modi2024joint} derives finite-sample error
rates for jointly learning related linear dynamical
systems. Although these contributions concern
different model classes and observation settings,
they provide useful foundations for understanding
estimation accuracy, regularization, and
identifiability in dynamical systems.

In parallel, stochastic optimization methods provide
computationally efficient alternatives for large-scale
finite-sum problems. Proximal stochastic gradient
methods with progressive variance reduction are
investigated in \citet{xiao2014proximal}, while
\citet{reddi2016proxsvrg} studies proximal
variance-reduced gradient methods for nonconvex
finite-sum optimization. These methods combine
stochastic gradient evaluations with variance
reduction and proximal regularization, making them
attractive for parameter-estimation problems
involving large datasets and structured penalties.
Recent developments also investigate variance
reduction in constrained finite-sum optimization.
In particular, \citet{gharehmohammadlou2026variance} develops
variance-reduced first-order methods for convex
finite-sum problems with nonlinear function
constraints and establishes gradient-complexity
guarantees. Although this framework addresses
a different optimization class, it demonstrates
the potential of variance-reduction techniques
for improving the computational efficiency of
finite-sum optimization.

Motivated by these developments, the present work
considers switched-system identification with
measurement errors in both inputs and outputs,
latent mode assignments, and unknown mode count
and subsystem orders. We formulate an
EIV identification problem that
combines a variance-normalized fitting criterion,
temporally regularized soft assignments,
mode-sparsity regularization, and group-sparse
lag selection. The resulting nonconvex problem
is addressed by an alternating optimization
procedure involving convex mode-assignment
updates, mode-weight optimization, and finite-sum
proximal stochastic variance-reduced gradient
(Prox-SVRG) steps.

Unlike approaches that explicitly enumerate all
candidate system structures, the proposed formulation
aims to identify the relevant modes and lag orders
within an initially model.
Its performance is investigated through numerical
comparisons with generalized total least squares
and dynamic programming (GTLS--DP), focusing on
structural recovery, mode-classification accuracy,
parameter-estimation error, and computational
scalability under increasing structural uncertainty.

\label{sec:notation}

\noindent
{\textbf{\textit{Notation.}}
Let $\mathbb{R}_{\geq 0}$ denote the nonnegative real numbers.
All vectors are column vectors; $(\cdot)^\top$ denotes transpose,
and $\|\cdot\|_2$ and $\|\cdot\|_F$ denote the Euclidean and
Frobenius norms, respectively. We use $I_d$,
$\operatorname{diag}(\cdot)$, $\mathbf{1}_{\{\cdot\}}$,
$|\mathcal S|$, $(a)_+$, and $e_k$ for the identity matrix,
(block-)diagonal operator, indicator function, set cardinality,
positive part of $a$, and the $k$th standard basis vector,
respectively. The symbol $N$ denotes the number of observations,
$t$ indexes discrete time, and $k$ indexes candidate operating
modes. The quantities $(\bar K,\bar n_a,\bar n_b)$ represent
assumed upper limits on the number of modes, output order,
and input order, respectively. Other model-specific notation,
including $\Theta$, $Q$, $\pi$, and the proximal operator,
is defined where it is introduced in the following sections.

\section{EIV Switch-System Identification}
\label{sec:eiv_problem}
\label{sec:problem}
We consider a switched autoregressive model with
exogenous inputs (SARX), whose dynamics are described by
\begin{equation}\!\!\!\!\!
y_t =
\sum_{j=1}^{n_a} a_{s_t,j} y_{t-j}
\!+\!
\sum_{j=1}^{n_b} b_{s_t,j} u_{t-j},
\quad s_t\in\{1,\ldots,K^\circ\},\!
\label{eq:true_sarx}
\end{equation}
where $u_t,y_t\in\mathbb{R}$ denote the input and output
signals at time $t$, respectively, and $s_t$ indicates
the active operating mode. Moreover, $K^\circ$ denotes
the number of operating modes, while $n_a$ and $n_b$
are the output and input orders, respectively.
For each mode $k=1,\ldots,K^\circ$, the coefficients
$\{a_{k,j}\}_{j=1}^{n_a}$ and
$\{b_{k,j}\}_{j=1}^{n_b}$ characterize the corresponding
subsystem, whose parameters are collected in
\begin{equation}
\theta_k =
[a_{k,1},\ldots,a_{k,n_a},
b_{k,1},\ldots,b_{k,n_b}]^\top.
\end{equation}
The switching sequence $S=(s_t)$, the number of
operating modes $K^\circ$, and the orders $(n_a,n_b)$
are assumed to be \emph{unknown}.

We consider an errors-in-variables (EIV) setting,
where the input and output signals are observed
through noisy measurements as
\begin{equation}
\left\{
\begin{aligned}
\widetilde y_t &= y_t+v_t,\\
\widetilde u_t &= u_t+w_t,
\end{aligned}
\right.
\qquad \text{for } t=0,\ldots,N-1,
\label{eq:eiv_measurements}
\end{equation}
where $N$ denotes the
number of observations. The measurement errors
$v_t$ and $w_t$ are assumed to be mutually
independent, zero-mean, and independent across
time, with known variances $\sigma_y^2$ and
$\sigma_u^2$, respectively.
Furthermore, we assume that structural upper
bounds $(\bar K,\bar n_a,\bar n_b)$ are available,
satisfying
$K^\circ\leq\bar K$, $n_a\leq\bar n_a$, and
$n_b\leq\bar n_b$.

\noindent\textbf{\textit{Identification Problem.}}
Given the noisy observations
$\mathcal D=\{(\widetilde u_t,\widetilde y_t)\}_{t=0}^{N-1}$,
the noise variances $\sigma_y^2,\sigma_u^2$, and the
structural upper bounds $(\bar K,\bar n_a,\bar n_b)$,
estimate the number of operating modes $K^\circ$,
the subsystem orders $(n_a,n_b)$, the subsystem
parameters $\{\theta_k\}_{k=1}^{K^\circ}$,
and the latent switching sequence $S=(s_t)$,
up to a permutation of mode labels.

The proposed identification scheme is presented in
Section~\ref{sec:identification_scheme}.

\section{Identification Models and Algorithms}
\label{sec:identification_scheme}

We develop an EIV switch-system identification model and an
alternating optimization scheme. Each component of the objective is
introduced separately to clarify its role. We then explain the parameter
and mode updates and relate the resulting method to an exhaustive
GTLS--DP benchmark.

\subsection{EIV switch-system Model}
\label{subsec:identification_models}
\label{subsec:eiv_scheme}

We use $\bar K$ candidate modes and maximum lag orders
$(\bar n_a,\bar n_b)$, assumed to contain the true structure.
Let vector $\Theta$ be defined as 
\begin{equation}
 \Theta=[\theta_1,\ldots,\theta_{\bar K}]^\top,   
\end{equation}
with
$ 
\theta_k=[\theta_{k,y}^\top,\theta_{k,u}^\top]^\top.    
$ 
For any $t$, let also $\widetilde\varphi_t$ be defined as the vector of past $\bar{n}_a$ measurements of output and past $\bar{n}_b$ measurements of input, i.e., we have  
\begin{equation}
\widetilde\varphi_t=
[\widetilde y_{t-1},\ldots,\widetilde y_{t-\bar n_a},
 \widetilde u_{t-1},\ldots,\widetilde u_{t-\bar n_b}]^\top.
\label{eq:eiv_regressor}
\end{equation}
We additionally define $t_0$ as $t_0=\max\{\bar n_a,\bar n_b\}$, which indicates the number of lags considered in the model.

\noindent\textbf{\textit{EIV Fitting.}}
Ordinary least squares treats the regressor as noise-free, which is
inappropriate when delayed outputs and inputs are measured with error.
For a sample generated by mode $k$, substituting the measurement model
\eqref{eq:eiv_measurements} into \eqref{eq:true_sarx} results in
\begin{equation}
 r_{tk}=v_t-\sum_{j=1}^{n_a}a_{k,j} v_{t-j}
             -\sum_{j=1}^{n_b}b_{k,j} w_{t-j}.
\label{eq:eiv_residual}
\end{equation}
Under the independent, zero-mean measurement-noise assumption, we have
\begin{equation}
\operatorname{Var}(r_{tk})=
\sigma_y^2\bigl(1+\|a_k\|_2^2\bigr)
+\sigma_u^2\|b_k\|_2^2.
\label{eq:eiv_variance}
\end{equation}
In \eqref{eq:eiv_variance}, the first term accounts for the current and delayed output errors, and the second term accounts for the errors in delayed inputs. Replacing the
unknown true coefficients by candidate coefficients motivates introducing a loss function defined as
\begin{equation}
\ell_{tk}(\theta_k)=\frac{1}{2}
\frac{\|\widetilde y_t-\widetilde\varphi_t^\top\theta_k\|_2^2}
{\sigma_y^2(1+\|\theta_{k,y}\|_2^2)
 +\sigma_u^2\|\theta_{k,u}\|_2^2}.
\label{eq:eiv_loss}
\end{equation}
This is a per-sample variance-normalized EIV criterion, not the
full joint likelihood: neighboring Switch-system residuals can share delayed
measurement errors. It gives the data-fitting component
\begin{equation}
J_{\mathrm{fit}}=\frac{1}{n}
\sum_{t\in\mathcal I}\sum_{k=1}^{\bar K}
q_{tk}\ell_{tk}(\theta_k),
\label{eq:fit_term}
\end{equation}
where $\mathcal I=\{t_0,\ldots,N-1\}$ and $n=|\mathcal I|$.

\noindent\textbf{\textit{Soft Assignments \& Temporal Consistency.}}
Let $q_t=(q_{t1},\ldots,q_{t\bar K})^\top
\in\Delta_{\bar K}$ represent the membership weights of
observation $t$, and let $\pi\in\Delta_{\bar K}$ denote
the global candidate-mode weights, where
$\Delta_{\bar K}=\{x\geq0:\sum_k x_k=1\}$.
We introduce the relative-entropy penalty
\begin{equation}
J_{\mathrm{ent}}=\frac{T}{n}
\sum_{t\in\mathcal I}\sum_{k=1}^{\bar K}
q_{tk}\log\frac{q_{tk}}{\pi_k},\qquad T>0.
\label{eq:entropy_term}
\end{equation}
This term discourages premature hard assignments while
linking local memberships to global mode occupancy.
As $T$ decreases, assignments become more selective.
We adopt $0\log 0=0$ and $q_{tk}=0$ when $\pi_k=0$.
However, relative entropy alone does not account for
the temporal persistence of switched systems, whose
operating modes typically remain unchanged over
consecutive observations. We therefore introduce
a temporal regularization term that penalizes abrupt
changes between successive assignment vectors, defined as
\begin{equation}
J_{\mathrm{time}}=
\frac{\gamma}{2(n-1)}
\sum_{t>t_0}\|q_t-q_{t-1}\|_2^2,
\qquad \gamma\geq 0.
\label{eq:temporal_term}
\end{equation}

\noindent\textbf{\textit{Mode Sparsity.}}
An oversized candidate set should not force all $\bar K$ modes
to remain active. We therefore use
\begin{equation}
J_{\mathrm{mode}}=\lambda_K
\sum_{k=1}^{\bar K}\sqrt{\pi_k+\varepsilon},
\qquad\lambda_K\geq0,\ \varepsilon>0.
\label{eq:mode_penalty}
\end{equation}
Since the square root is concave, minimizing this term on the
simplex favors concentrated mode weights and suppresses
unnecessary modes.

\noindent\textbf{\textit{Lag-group Sparsity.}}
To select switch-system orders, collect the coefficients of output lag $j$
across all candidate modes in $\Theta_{a,j}$, and similarly
define $\Theta_{b,j}$. We impose
\begin{equation}
\begin{aligned}
R(\Theta)={}&\sum_{j=1}^{\bar n_a}
 \lambda_{a,j}\|\Theta_{a,j}\|_2+\sum_{j=1}^{\bar n_b}
 \lambda_{b,j}\|\Theta_{b,j}\|_2.
\end{aligned}
\label{eq:lag_penalty}
\end{equation}
Each penalty acts on a complete lag across modes, so that an
irrelevant lag can be removed jointly rather than coefficient by
coefficient. In the structural-scaling experiment, the adopted
weights are $\lambda_{\cdot,1}=0$, $\lambda_{\cdot,2}=0.002$,
and $\lambda_{\cdot,j}=0.12$ for $j\geq3$; these encode a
strong low-order prior. 

\noindent\textbf{\textit{Complete Identification Criterion.}}
Combining the preceding components gives
\begin{equation}\!\!\!\!\!\!\!\!\!\!
\begin{aligned}
\min_{\Theta,Q,\pi}\quad
J(\Theta,Q,\pi)={}&J_{\mathrm{fit}}+J_{\mathrm{ent}}
 +J_{\mathrm{time}}\\
 &+J_{\mathrm{mode}}+R(\Theta)\\
\text{s.t.}\quad &q_t\in\Delta_{\bar K}, \qquad \forall t\in\mathcal I,\\
 &\pi\in\Delta_{\bar K}.
\end{aligned}
\label{eq:proposed_problem}
\end{equation}
The terms respectively promote EIV-consistent data fitting,
assignment softness, temporal persistence, a small number of
active modes, and parsimonious orders. The complete objective
is nonconvex because it couples assignments and parameters,
contains a coefficient-dependent EIV denominator, and includes
the concave mode-sparsity penalty. We consequently seek a
stationary solution by alternating block updates rather than
claiming global optimality. After fitting, define
\begin{equation}
\begin{aligned}
\widehat K&=\#\{k:\pi_k>\tau_K\},\\
\widehat n_a&=\max\{j:\|\Theta_{a,j}\|_2>\tau_p\},\\
\widehat n_b&=\max\{j:\|\Theta_{b,j}\|_2>\tau_p\}.
\end{aligned}
\label{eq:structure_selection}
\end{equation}
Here an empty maximum is interpreted as zero; the numerical
scaling experiment uses $\tau_K=0.02$ and $\tau_p=0.05$.

\subsection{Finite-Sum Prox-SVRG Identification Algorithm}
\label{subsec:identification_algorithms}
\label{subsec:svrg}

We alternate updates of $Q$, $\pi$, and $\Theta$ while retaining
the other two blocks fixed. The corresponding subproblems and
computations are detailed below.

\noindent\textit{The assignment ($Q$) update.}

For fixed $(\Theta,\pi)$, the terms depending on $Q$ yield
\begin{equation}
\begin{aligned}
Q^+\in\argmin_{q_t\in\Delta_{\bar K}}\quad
&\frac1n\sum_{t,k}q_{tk}\ell_{tk}(\theta_k)
+\frac{T}{n}\sum_{t,k}q_{tk}\log\frac{q_{tk}}{\pi_k}
\\&\qquad+\frac{\gamma}{2(n-1)}
  \sum_{t>t_0}\|q_t-q_{t-1}\|_2^2.
\end{aligned}
\label{eq:q_update}
\end{equation}
For $T>0$ and positive $\pi$, this is a strictly convex
simplex-constrained problem. Its temporal coupling is handled
by alternating updates of even and odd time indices. With
neighboring assignment vectors fixed, a row $q_t$ solves
\begin{equation}
\begin{aligned}
\min_{q\in\Delta_{\bar K}}\quad
&\sum_kq_k\ell_{tk}+T\sum_kq_k\log\frac{q_k}{\pi_k}\\
&\qquad\qquad+\frac{\gamma n}{2(n-1)}
 \sum_{j\in\mathcal N_t}\|q-q_j\|_2^2,
\end{aligned}
\label{eq:q_row}
\end{equation}
where $\mathcal N_t$ contains the available adjacent indices.
Each row is solved numerically under the simplex constraint;
when $\gamma=0$, it reduces to
\begin{equation}
q_{tk}=\frac{\pi_k\exp(-\ell_{tk}/T)}
{\sum_h\pi_h\exp(-\ell_{th}/T)}.
\label{eq:q_softmax}
\end{equation}
The latter expression illustrates the effect of the EIV cost
and assignment temperature, but the experiments retain the
temporal penalty $\gamma\geq0$.

\noindent\textbf{\textit{Mode-weight Update.}}
Let $\bar q_k=n^{-1}\sum_tq_{tk}$ be the mean membership of
mode $k$. With $Q$ fixed, the relevant objective is
\begin{equation}\!\!\!\!\!
\pi^+\in\argminOp_{\pi\in\Delta_{\bar K}}
\left\{-T\sum_k\bar q_k\log\pi_k
+\lambda_K\sum_k\sqrt{\pi_k+\varepsilon}\right\}\!.
\label{eq:pi_update}\!\!\!\!\!
\end{equation}
The entropy contribution discourages eliminating populated
modes, whereas the concave penalty favors a sparse mode set.
The subproblem is generally nonconvex; we optimize a
softmax parameterization of $\pi$ numerically, obtaining a
local solution, and then update $Q$ again by \eqref{eq:q_update}.

\noindent\textbf{\textit{Parameter Update.}}
For fixed $(Q,\pi)$, the parameter objective has finite-sum
composite form
\begin{equation}
\min_{\Theta}\ F(\Theta)+R(\Theta),
\label{eq:theta_update}
\end{equation}
where we have 
\begin{equation}
F(\Theta)=\frac1n\sum_{t\in\mathcal I}f_t(\Theta),
\end{equation}
and $f_t(\Theta)=\sum_kq_{tk}\ell_{tk}(\theta_k)$.
Let $e_{tk}=\widetilde y_t-\widetilde\varphi_t^\top\theta_k$,
$d_k=\sigma_y^2(1+\|\theta_{k,y}\|^2)
+\sigma_u^2\|\theta_{k,u}\|^2$, and
$D=\operatorname{diag}(\sigma_y^2I_{\bar n_a},
\sigma_u^2I_{\bar n_b})$. The gradient used in the update is
\begin{equation}
\nabla_{\theta_k}\ell_{tk}
=-\frac{e_{tk}\widetilde\varphi_t}{d_k}
 -\frac{e_{tk}^2 D\theta_k}{d_k^2}.
\label{eq:eiv_gradient}
\end{equation}
At each epoch, form a snapshot $\widetilde\Theta$ and its
full gradient $\nabla F(\widetilde\Theta)$. For mini-batch
$\mathcal B_m$, define
\begin{equation}
\begin{aligned}
v_m={}&\nabla F_{\mathcal B_m}(\Theta^m)
-\nabla F_{\mathcal B_m}(\widetilde\Theta)
+\nabla F(\widetilde\Theta),
\end{aligned}
\label{eq:svrg_gradient}
\end{equation}
where $F_{\mathcal B}$ is the average of $f_t$ over the batch.
The Prox-SVRG parameter step is
\begin{equation}
\Theta^{m+1}=\operatorname{prox}_{\eta, R}
  (\Theta^m-\eta v_m).
\label{eq:svrg_prox}
\end{equation}
For any output- or input-lag group $g$ with penalty $\lambda_g$,
the proximal operator in \eqref{eq:svrg_prox} applies the
group soft-thresholding rule
\begin{equation}
\operatorname{prox}_{\eta,\lambda_g}(z_g)
=\left(1-\frac{\eta\lambda_g}{\|z_g\|_2}\right)_+z_g,
\label{eq:group_prox}
\end{equation}
with value zero when $z_g=0$. Hence, stochastic variance
reduction addresses the smooth EIV term while the proximal
step removes irrelevant lag groups. The assignments are then
refreshed using \eqref{eq:q_update}; the temperature $T$ is
reduced between epochs to sharpen mode separation.

\noindent\textbf{\textit{From soft to hard mode assignments.}}
In the fitting jump models framework \citep{bemporad2018jump},
subsystem estimation alternates with optimization of a discrete
switching sequence. For one-hot $q_t=e_{s_t}$, the temporal term
\eqref{eq:temporal_term} satisfies
\begin{equation}
\frac{\gamma}{2(n-1)}\sum_{t>t_0}
\|q_t-q_{t-1}\|_2^2
=\frac{\gamma}{n-1}\sum_{t>t_0}
\mathbf1_{\{s_t\ne s_{t-1}\}}.
\label{eq:fjm_connection}
\end{equation}
Thus, our temporal regularizer continuously extends a jump
penalty. A hard-assignment EIV counterpart is obtained by
alternating GTLS subsystem fits and DP mode assignments.
The benchmark repeats this process over candidate structures,
whereas the proposed method estimates structural sparsity
within one model.
\subsection{Generalized Total Least Squares
and Dynamic Programming}
\label{subsec:gtls}

\noindent\textbf{\textit{GTLS fitting term.}}
For a candidate $(K,n_a,n_b)$, define
\begin{equation}
\begin{gathered}
z_t=[\widetilde y_t,\widetilde y_{t-1},\ldots,
\widetilde y_{t-n_a},\widetilde u_{t-1},\ldots,
\widetilde u_{t-n_b}]^\top,\\
h_k=[1,-\theta_k^\top]^\top,\quad
\Sigma=\operatorname{diag}(\sigma_y^2I_{n_a+1},
\sigma_u^2I_{n_b}).
\end{gathered}
\label{eq:gtls_definition}
\end{equation}
Since $z_t^\top h_k$ is the measured regression residual,
its variance-normalized squared cost is
\begin{equation}
c_{tk}=\frac{(z_t^\top h_k)^2}{h_k^\top\Sigma h_k}.
\label{eq:gtls_cost}
\end{equation}
For fixed hard assignments $S$, let $\mathcal I_k$ be the
indices assigned to mode $k$ and set
$A_k=|\mathcal I_k|^{-1}\sum_{t\in\mathcal I_k}z_tz_t^\top$.
The GTLS coefficient update is given by the generalized
smallest-eigenvalue problem
\begin{equation}
A_kh_k=\lambda_k\Sigma h_k,\quad
\lambda_k=\lambda_{\min}(A_k,\Sigma),
\label{eq:gtls_eigen}
\end{equation}
followed by scaling $h_{k,1}=1$, provided its first component
is nonzero, and recovering the subsystem coefficients.

\noindent\textbf{\textit{The dynamic-programming assignment.}}
With the GTLS coefficients fixed, the mode sequence minimizes
\begin{equation}
\min_{S}\left\{\sum_{t}c_{t,s_t}
+\gamma_{\mathrm{DP}}\sum_{t>t_0}
\mathbf1_{\{s_t\ne s_{t-1}\}}\right\}.
\label{eq:gtls_dp}
\end{equation}
Let $V_t(k)$ denote the minimum cumulative cost of a
sequence ending in mode $k$ at time $t$. The DP recurrence is
\begin{equation}\left\{
\begin{aligned}
&V_{t_0}(k)\,=\,c_{t_0k},\\
&V_t(k)\,=\,c_{tk}+\min_h\left\{V_{t-1}(h)
+\gamma_{\mathrm{DP}}\mathbf1_{\{h\ne k\}}\right\}.
\end{aligned}\right.
\label{eq:dp_recursion}
\end{equation}
Backtracking from $\argmin_kV_{N-1}(k)$ recovers the
optimal mode sequence for the current coefficients. GTLS
\eqref{eq:gtls_eigen} and DP \eqref{eq:dp_recursion} are
alternated until assignment stabilization or a preset limit.

\noindent\textbf{\textit{Structural selection and complexity.}}
Each candidate is accepted if its mean normalized residual
$\overline c_{K,n_a,n_b}$ does not exceed
$\tau_{\mathrm{GTLS}}=1.05$; among accepted candidates, the
benchmark selects the smallest $K(n_a+n_b)$.
There are
\begin{equation}
M_{\mathrm{cand}}=\bar K\bar n_a\bar n_b
\label{eq:number_candidates}
\end{equation}
possible structures in this search. The benchmark therefore
solves a separate GTLS--DP problem for each candidate,
while Prox-SVRG solves a single larger problem.

\noindent\textbf{\textit{Algorithms.}}
Algorithm \ref{alg:proposed} alternates between updating the
soft mode assignments, global mode weights, and subsystem
parameters. The assignments are obtained by solving a
temporally regularized convex problem, while the mode
weights are updated to promote sparsity. For parameter
estimation, Prox-SVRG combines variance-reduced stochastic
gradients with proximal updates to suppress irrelevant
lag groups. The active modes and subsystem orders
are determined by thresholding the estimated mode weights
and lag-group norms. The GTLS--DP benchmark \ref{alg:gtls} exhaustively evaluates candidate structures defined by the admissible mode counts and
input-output orders. For each candidate, it alternates
between GTLS-based subsystem parameter estimation and
dynamic programming for optimal hard mode assignments,
using multiple initializations. The best solution is
retained for each structure, and the final model is
selected by minimizing structural complexity among
candidates satisfying a prescribed fitting threshold.

\begin{algorithm}[t]
\caption{Finite-Sum Prox-SVRG for Unknown-Structure
EIV-PWA Identification}
\label{alg:proposed}
\begin{algorithmic}[1]
\State Initialize $\Theta$, $\pi_k=1/\bar K$, $Q$,
and $T=T_0$
\For{$r=1,\ldots,R$}
  \State Evaluate $\ell_{tk}$ using \eqref{eq:eiv_loss}
  \State Update $Q$ by
  \eqref{eq:q_update}--\eqref{eq:q_row}
  \State Update $\pi$ by \eqref{eq:pi_update}
  \State Refresh $Q$ using \eqref{eq:q_update}
  \State Set $\widetilde{\Theta}\gets\Theta$ and compute
  $\widetilde g=\nabla F(\widetilde{\Theta})$
  \State Set $\Theta_0\gets\widetilde{\Theta}$
  \For{$m=0,\ldots,M-1$}
    \State Sample mini-batch $\mathcal B_m$
    \State Compute $v_m$ by \eqref{eq:svrg_gradient}
    \State Update $\Theta_{m+1}$ via the proximal step
    \eqref{eq:svrg_prox}--\eqref{eq:group_prox}
  \EndFor
  \State Set $\Theta\gets\Theta_M$
  \State Recompute $\ell_{tk}$ and refresh $Q$
  using \eqref{eq:q_update}
  \State Set $T\gets\max\{T_{\min},\rho T\}$
\EndFor
\State Select the estimated structure using
\eqref{eq:structure_selection}
\end{algorithmic}
\end{algorithm}

\begin{algorithm}[t]
\caption{Exhaustive GTLS--DP Structure Selection}
\label{alg:gtls}
\begin{algorithmic}[1]
\For{$K=1,\ldots,\bar K$}
  \For{$n_a=1,\ldots,\bar n_a$}
    \For{$n_b=1,\ldots,\bar n_b$}
      \For{$j=1,\ldots,J_{\mathrm{init}}$}
        \State Initialize $K$ homogeneous normals $h_k$
        \Repeat
          \State Compute costs using \eqref{eq:gtls_cost}
          \State Update assignments using \eqref{eq:dp_recursion}
          \State Update normals using \eqref{eq:gtls_eigen}
        \Until{assignments converge}
      \EndFor
      \State Keep the best GTLS--DP solution
      \State Compute mean residual $\overline c_{K,n_a,n_b}$
    \EndFor
  \EndFor
\EndFor
\State Define the admissible structure set
\Statex \hspace{\algorithmicindent}
$\mathcal A=\{(K,n_a,n_b):
\overline c_{K,n_a,n_b}\leq\tau_{\mathrm{GTLS}}\}$
\If{$\mathcal A\neq\emptyset$}
  \State Select $(\widehat K,\widehat n_a,\widehat n_b)$
  from $\mathcal A$ by lexicographically minimizing
  \Statex \hspace{\algorithmicindent}
  $\bigl(K(n_a+n_b),K,\overline c_{K,n_a,n_b}\bigr)$
\Else
  \State Report no admissible structure
\EndIf
\end{algorithmic}
\end{algorithm}

\section{Numerical Results}
\label{sec:numerical_results}
We compare finite-sum Prox-SVRG with exhaustive GTLS--DP
on a $2\times2$  EIV switch system with true structure
$(K^\circ,n_a^\circ,n_b^\circ)=(3,2,2)$.
We consider 17 structural envelopes with
$M=\bar K\bar n_a\bar n_b$ candidate structures,
ranging from 45 to 686.
We use 30 Monte Carlo realizations of length $N=4800$,
with fixed subsystem matrices and varying inputs, switching
sequences, and Gaussian measurement errors whose standard
deviations equal $20\%$ of the corresponding signal standard
deviations. Each realization is shared across both methods
and all envelopes. Prox-SVRG uses 60 epochs, 50 inner steps
per epoch, mini-batches of 96, and three restarts; GTLS--DP
uses four restarts. We report runtime, mode-classification
accuracy, and mean Frobenius parameter error after optimal
mode matching and zero-padding of estimated lag matrices.

GTLS--DP has lower mean runtime through $M=546$, while
Prox-SVRG reduces mean runtime by $7.44\%$ at $M=637$ and
$8.83\%$ at $M=686$ (Fig.~\ref{fig:eiv_runtime}).
At $M=686$, the respective mean runtimes are $214.95$
and $235.77$ s. Across 510 paired comparisons, exact
structural recovery is $508/510$ for Prox-SVRG and
$510/510$ for GTLS--DP; the two Prox-SVRG failures
occur at $M=504$ and $546$.

\begin{figure}[t]
\centering
\includegraphics[width=\columnwidth]{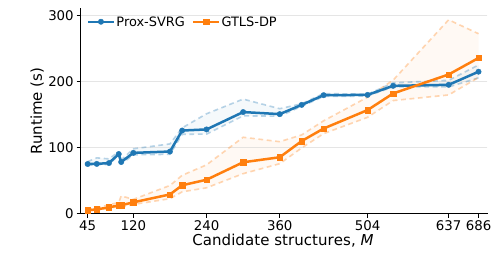}
\label{fig:eiv_runtime}
\end{figure}

GTLS--DP attains slightly higher mean mode-classification
accuracy: $99.853\%$ versus $99.788\%$ for Prox-SVRG
(Fig.~\ref{fig:eiv_mode_accuracy}). Its mean parameter
error is also lower: $0.04965$ versus $0.05499$
(Fig.~\ref{fig:eiv_parameter_error}). Thus, the runtime
advantage at the two largest envelopes comes with a
small loss in estimation accuracy.

\begin{figure}[t]
\centering
\includegraphics[width=\columnwidth]{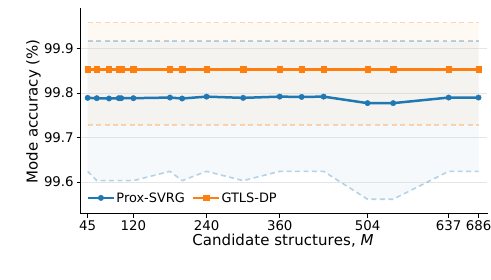}
\caption{Mean mode-classification accuracy over 30 repetitions (zoomed scale).}
\label{fig:eiv_mode_accuracy}
\end{figure}

\begin{figure}[t]
\centering
\includegraphics[width=\columnwidth]{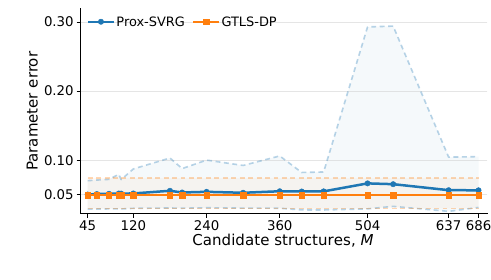}
\caption{Mean Frobenius parameter error after optimal mode matching.}
\label{fig:eiv_parameter_error}
\end{figure}

Because the envelopes are not all nested, these comparisons
describe empirical scaling rather than isolate the effect
of a single structural dimension.

\section{Concluding remarks and future work}
\label{sec:conclusion}
This paper proposes an errors-in-variables (EIV) identification framework for switched dynamical systems with unknown mode count and subsystem orders. The proposed approach combines temporally regularized soft mode assignments, sparsity-based structure selection, and alternating optimization incorporating proximal stochastic variance-reduced gradient updates. Numerical comparisons across 17 structural envelopes and 30 Monte Carlo repetitions demonstrate high structural-recovery accuracy for both methods, with recovery rates of 99.61 for Prox-SVRG and 100 for exhaustive GTLS--DP. Although GTLS--DP achieves slightly higher mode-classification accuracy and lower parameter-estimation error, Prox-SVRG reduces mean runtime by 7.44 and 8.83 at the two largest structural envelopes, respectively. These results highlight the favorable computational efficiency of the proposed approach under increasing structural uncertainty, with a modest trade-off in estimation accuracy.

\bibliography{ifacconf_noDOI}
\end{document}